\documentclass[intlimits]{article}
\usepackage{latexsym,amsfonts,amssymb,amsmath,amsthm,amscd}

\theoremstyle{plain}

\newcommand\RR{\ensuremath{\mathbb{R}}}

\DeclareMathOperator{\trace}{tr}

\begin{document}
\title{A simple proof of monotonicity for linear cooperative systems
of ODEs}

\author{Janusz Mierczy\'nski}


\maketitle

\begin{abstract}
We present a simple proof of monotonicity for cooperative systems of
linear ordinary differential equations, without having recourse to
approximation by strongly cooperative systems.
\end{abstract}

A system of linear ordinary differential equations (ODEs)
\begin{equation}
\label{ODE}
x' = A(t) x, \qquad t \in (a, b),
\end{equation}
where $A(t) = [a_{ij}(t)]_{i,j=1}^n$ is a matrix with real entries,
is called {\em cooperative\/} if for each $t \in (a, b)$ the
off-diagonal entries of the matrix $A(t)$ are nonnegative.

We assume in the remainder that the matrix function $A \colon (a, b)
\to \RR^{n \times n}$ is continuous, so for any $(t_0, x_0) \in (a,
b) \times \RR^n$ there exists a unique solution
\begin{equation*}
[\, (a, b) \ni t \mapsto x(t; t_0, x_0) \in \RR^n \, ]
\end{equation*}
of the initial value problem
\begin{equation}
\label{IVP}
\begin{cases}
x' = A(t) x, \\
x(t_0) = x_0.
\end{cases}
\end{equation}

An important property of cooperative systems of ODEs is that their
solution operators have some monotonicity properties.  To formulate
them, write
\begin{equation*}
\begin{aligned}
\RR^{n}_{+} := {} & \{\, x = (x_1, \dots, x_n) \in \RR^{n} : x_i \ge
0 \text{ for all } 1 \le i \le n \, \}, \\
\RR^{n}_{++} := {} & \{\, x = (x_1, \dots, x_n) \in \RR^{n} : x_i > 0
\text{ for all } 1 \le i \le n \, \}.
\end{aligned}
\end{equation*}

\begin{enumerate}
\item[{\textbf{(M1)}}]
{\em For any $t_0 \in (a, b)$, if $x_0 \in \RR^n_{+}$ then $x(t;
t_0, x_0) \in \RR^n_{+}$ for all $t \in (t_0, b)$.}
\item[{\textbf{(M2)}}]
{\em For any $t_0 \in (a, b)$, if $x_0 \in \RR^n_{++}$ then $x(t;
t_0, x_0) \in \RR^n_{++}$ for all $t \in (t_0, b)$.}
\end{enumerate}

The approach to proving the monotonicity (more precisely, property
\textbf{(M1)}) in all the literature on cooperative systems I have
consulted (cf., e.g., \cite{Coppel}, \cite{L-L}, \cite{Hirsch-Smith},
\cite{Walter}) rests on the following.
\begin{itemize}
\item
A cooperative system is approximated by systems whose matrices
have their off-diagonal entries positive ({\em strongly
cooperative systems\/}).
\item
It is proved that in a strongly cooperative system, if a nonzero
$x_0$ is in $\RR^n_{+}$ then $x(t; t_0, x_0) \in \RR^n_{++}$ for
$t \in (t_0, b)$. That property is called {\em strong
monotonicity\/}.
\item
With the help of the continuous dependence of a solution on a
parameter one obtains \textbf{(M1)} for the original system.
\end{itemize}
Indeed, the above approach is usually applied to general nonlinear
systems, and from that point of view linear systems are special cases
of nonlinear ones.

When one wants to obtain (the apparently stronger) property
\textbf{(M2)} one needs to apply some additional ideas, for instance,
that the solution mapping is a homeomorphism (for a different method,
see \cite{Coppel}).

\bigskip
Below I present a simple and elementary approach, to the best of my
knowledge (strangely) overlooked so~far.

\begin{proof}[Proof of \textbf{(M2)}]
Fix $t_0 \in (a, b)$, $x_0 \in \RR^n_{++}$, and write $x(t) =
(x_1(t), \dots, x_n(t))$ for $x(t; t_0, x_0)$.  We have
\begin{equation*}
\begin{aligned}
\frac{d}{dt} \prod_{i} x_i = {} & \sum_{i} \Bigl( x'_{i}  \prod_{j
\ne i} x_j  \Bigr) = \sum_{i} \biggl( \Bigl( \sum_{k} a_{ik}(t) x_{k}
\Bigr)
\prod_{j \ne i} x_j  \biggr) \\
\ge {} & \sum_{i} \Bigl( a_{ii}(t) \, x_{i} \prod_{j \ne i} x_j
\Bigr) = \trace{A(t)}  \prod_{i} x_i.
\end{aligned}
\end{equation*}
So, denoting $\xi(t) := x_{1}(t) \cdot \hdots \cdot x_{n}(t)$, there
holds
\begin{equation*}
\xi'(t) \ge \trace{A(t)} \, \xi(t), \qquad t \in (a, b),
\end{equation*}
as well as $\xi(t_0) > 0$.  By continuity, $\xi(t) > 0$ for $t > t_0$
sufficiently close to $t_0$.  Suppose to the contrary that there are
$\theta \in (t_0, b)$ at which $\xi(\theta) = 0$.  Let $\tau \in
(t_0, b)$ be the least such $\theta$.  Since, for each $t \in (t_0,
\tau)$ there holds
\begin{equation*}
\xi(t) \ge \xi(t_0) \exp\Bigl(\int\limits_{t_0}^{t} \trace{A(s)} \,
ds\Bigr),
\end{equation*}
by letting $t \nearrow \tau$ we obtain a contradiction, as the
left-hand side tends to zero whereas the right-hand side goes to a
positive number.

It suffices now to notice that this means that $x_i(t) > 0$ for all
$t \in (t_0, b)$ and all indices $i \in \{1, \dots, n\}$.
\end{proof}

It should be mentioned that I came across the idea of using the
product of the coordinates while reading the
paper~\cite{Benaim-Schreiber} by M. Bena\"{\i}m and S. J. Schreiber.

\begin{enumerate}
\item
As the mathematical tools used are rather elementary, the proof
is accessible to students who have taken a course in elementary
calculus, linear algebra and differential equations.

\item
The proof carries over, almost word by word, to the case of
Carath\'eodory conditions.

\item
It would be interesting to give an analog of it for the case of
systems of linear ODEs preserving an arbitrary cone in $\RR^{n}$,
compare e.g. \cite{Walcher}.
\end{enumerate}

\end{document}